\documentclass[10pt]{article}
\usepackage{amssymb}
\usepackage[fleqn]{amsmath}
\usepackage{mathrsfs}
\usepackage{amsfonts}
\usepackage{tikz}
\usepackage[numbers,sort&compress]{natbib}

\usepackage{amsmath,amsthm,amsfonts,amssymb,amscd}

\allowdisplaybreaks[4]
\newtheorem {Problem} {Problem}[section]
\newtheorem {Theorem} [Problem]{Theorem}
\newtheorem {Lemma}[Problem]{Lemma}
\newtheorem{Conjecture}[Problem]{Conjecture}
\newenvironment {Proof}{\noindent {\bf Proof.}}{\hfill\ensuremath{\square}}
\newcommand*{\QEDB}{\hfill\ensuremath{\square}}

\begin{document}

\title{The signless Laplacian spectral radius of graphs  without intersecting odd cycles   \thanks{This work is supported by  the National Natural Science Foundation of China (Nos. 11971311, 12026230),
% the Montenegrin-Chinese Science and Technology Cooperation Project (No.3-12),
 and the Hainan Provincial Natural Science Foundation of China (Nos. 120RC453, 120MS002).
\newline \indent Email: mzchen@hainanu.edu.cn, amliu@hainanu.edu.cn, xiaodong@sjtu.edu.cn. \newline  \indent $^{\dagger}$Corresponding author:
Xiao-Dong Zhang (Email: xiaodong@sjtu.edu.cn),}}

\author{ Ming-Zhu Chen, A-Ming Liu,\\
School of Science, Hainan University, Haikou 570228, P. R. China, \\
\and  Xiao-Dong Zhang$^{\dagger}$%\footnote{Corresponding author. E-mail: xiaodong@sjtu.edu.cn}
\\School of Mathematical Sciences, MOE-LSC, SHL-MAC\\
Shanghai Jiao Tong University,
Shanghai 200240, P. R. China}

\date{}
\maketitle

\begin{abstract}
Let  $F_{a_1,\dots,a_k}$ be a graph consisting of $k$ cycles of odd length $2a_1+1,\dots, 2a_k+1$, respectively which intersect in exactly a common vertex, where $k\geq1$ and $a_1\ge a_2\ge \cdots\ge a_k\ge 1$.  In this paper, we
present a sharp upper bound for the signless Laplacian spectral radius of all $F_{a_1,\dots,a_k}$-free graphs and characterize all extremal graphs which attain the bound. The stability methods and structure  of
 graphs associated with the eigenvalue are adapted for the proof.
\end{abstract}

\section{Introduction}
 Let $G$ be an undirected simple graph with vertex set
$V(G)=\{v_1,\dots,v_n\}$ and edge set $E(G)$, where $e(G)$ is the number of edges of $G$.
%Denote  $\delta(G)$ by   minimum degree of $G$.
For $v\in V(G)$,  the \emph{neighborhood} $N_G(v)$ of $v$  is $\{u: uv\in E(G)\}$ and the \emph{degree} $d_G(v)$ of $v$  is $|N_G(v)|$.
We write $N(v)$ and $d(v)$ for $N_G(v)$ and $d_G(v)$ respectively if there is no ambiguity.
 Denote  by $\triangle(G)$ and $\delta(G)$  the maximum and minimum degree of $G$, respectively. Denote by $P_n$ and $C_n$ the path and the cycle of order $n$, respectively.
%For a subgraph $H$  of $G$,  the \emph{$H$-neighborhood}  $N_H(v)$ of $v$  is $N_G(v)\bigcap V(H)$ and the
%\emph{$H$-degree} $d_H(v)$ of $v$  is $|N_H(v)|$.
%For two graphs $G$ and $H$, write $H\subseteq G$ if $G$ contains $H$ as a subgraph.
For $A, B\subseteq V(G)$,  $e(A)$ denotes the number of the  edges of  $G$ with both endvertices in $A$ and $e(A,B)$ denotes the number of the  edges of  $G$ with one endvertex in $A$ and the other  in $B$.
%In particular,  if $H$ is a proper subgraph of $G$, write $H\subseteq G$.
For two vertex disjoint graphs $G$ and $H$,  we denote by  $G\bigcup H$ and  $G\nabla H$  the \emph{union} of $G$ and $H$,
and the \emph{join} of $G$ and $H$, i.e., joining every vertex of $G$ to every vertex of $H$, respectively.
Denote by $k G$  the union of  $k$ disjoint  copies of $G$.
Denote by $S_{n,t}=K_t \bigcup \overline{K}_{n-t}$, i.e., the join of a complete graph $K_t$ and the independent set $I$ of size $n-t$.
Denote by  $L_{r,t}=K_1\nabla r K_t$, i.e., the graph consists of $r$ complete graph $K_{t+1}$ which intersect  in exactly a common vertex.
For graph notation and terminology undefined here, we refer the readers to \cite{BM}.

 We say a graph $G$ is \emph{$H$-free} if   it  does not contain  $H$ as a subgraph. The Tur\'{a}n number of a graph $H$ is the maximum number of edges that is in
an $H$-free graph of order $n$, and is denoted by $ex(n, H)$. An $H$-free
graph of order $n$  with $ex(n, H)$ edges is called an \emph{extremal
graph} for $H$ and  denote  $Ex(n, H)$ by  the set of all extremal graphs of order $n$ for $H$. To determine $ex(n, H)$ and characterize those graphs in $Ex(n, H)$ is a
fundamental problem (called Tur\'{a}n-type problem) in extremal graph theory.  It will be interesting to look for some nice graphs $H$  such that  $ex(n, H)$ and $Ex(n, H)$ will be characterized.  The graph  $F_{a_1,\dots,a_k}$ which  is a graph consisting of $k$ cycles of odd length $2a_1+1,\dots, 2a_k+1$, respectively which intersect in exactly a common vertex may be of interest, where $k\geq1$ and $a_1\ge a_2\ge\cdots\ge a_k\ge 1$. If $k=1$,  then $F_{a_1,\dots,a_k}$ is an odd cycle $C_{2a_1+1}$. Simonovits \cite{Simonovits1968} proved that  $ex(n, C_{2a_1+1})=\lfloor\frac{n^2}{4}\rfloor$ and $Ex (n, C_{2a_1+1})$ is a balance complete bipartite graph.
 % and $Ex(n, K_3)$ has been determined by Mantel \cite{}.
 If $k\ge 2$ and $a_1=\cdots=a_k=1$, then
 $F_{a_1,\dots,a_k}$ is denoted by $F^{(k)}$ which is called the \emph{friendship graph}.
 %, i.e., $S_{n,h}=K_h\nabla\overline{K}_{n-h}$.
  %, and
%reader is referred to \cite{} for surveys on this topic.
In  1995,  Erd\H{o}s,   F\"{u}redi,  Gould,  Gunderson \cite{Erdos1995}  significantly extended Mantel's result and proved the following interesting result.
\begin{Theorem}\label{thm1.1}\cite{Erdos1995}
Let  $k \geq 1$ and  $n \geq 50k^2$. Then
%
%if a graph $G$ of order $n$
%satisfies $e(G) > ex(n, F_k)$, then $G$ contains a copy of a $k$-fan, where
\[ ex(n,F^{(k)})=\Big\lfloor\frac{n^2}{4}\Big\rfloor+\left\{
\begin{array}{cccc}
 \vspace{1mm}
  k^2-k,&& \mbox{if $k$ is odd,}\\
   \vspace{1mm}
   k^2 -\frac{3}{2} k  ,&&\mbox{if $k$ is even.}
\end{array}\right.
\]
Furthermore, if $k$ is odd, then  $Ex(n, F^{(k)})$ consists of graphs which  are constructed by taking a complete bipartite graph with two partite of size $\lceil\frac{n}{2}\rceil$  and $\lfloor\frac{n}{2}\rfloor$ and embedding two vertex disjoint copies of $K_k$ in one side.
If $k$ is even, then  $Ex(n, F^{(k)})$ consists of graphs which are constructed by taking a complete bipartite graph with two partite of size $\lceil\frac{n}{2}\rceil$  and $\lfloor\frac{n}{2}\rfloor$ and embedding a graph with $2k-1$ vertices, $k^2 -\frac{3}{2} k $ edges with maximum degree $k-1$ in one partite.
\end{Theorem}
 If $k\ge 2$ and $a_1\ge \cdots \geq a_s\ge 2$, $a_{s+1}=\cdots=a_k=1$ with $1\le s\leq k$, then $F_{a_1,\cdots, a_k}$ is  denoted by  $H_{k,s}$,  i.e., $H_{k,s}$ is  the graph consisting of $k$ odd cycles and $k-s$ triangles which intersect in  exactly a common vertex.
 In 2018, Hou, Qiu and Liu \cite{Hou2018}, and Yuan \cite{Yuan2018} independently proved the following result.
 \begin{Theorem}\label{thm1.2}\cite{Hou2018,Yuan2018}
Let  $ k \geq 2$ and $1\le s\le k$.  Then
 $$ex(n,H_{k,s})=\Big\lfloor\frac{n^2}{4}\Big\rfloor+(k-1)^2$$
 for sufficient large $n$.
Moreover,  $Ex(n, H_{k,s})$ consists of a balance complete bipartite graph with a complete bipartite graph $K_{k-1,k-1}$ embedded in  one partite if $(k,s)\neq (4,1)$; a  balance complete bipartite graph with a complete bipartite graph $K_{3,3}$  or $3K_3$ embedded in one partite if $(k,s)=(4,1)$  \end{Theorem}

In spectral  extremal graph theory,  there is an  analogy  between the  Tur\'{a}n type problem  and the Brualdi-Solheid-Tur\'{a}n type problem which are proposed by  Nikiforov \cite{Nikiforov2010}. The Brualdi-Solheid-Tur\'{a}n type problem  is how to determine the maximum spectral radius of an $H$-free graph of order $n$ and characterize those graphs which attain the maximum spectral radius.
The Brualdi-Solheid-Tur\'{a}n type problem of the adjacent spectral spectral radius has been studied for various kinds of $H$ such as the complete
graph \cite{Nikiforov2007}, the complete bipartite graph \cite{Nikiforov2010-2},  the cycles or paths of specified
length \cite{Nikiforov2010}, the linear forest \cite{Chen2019-2}, and star forest \cite{Chen2021} and so on.  In addition, the Brualdi-Solheid-Tur\'{a}n type problem of the signless Laplacian spectral radius has also been investigated extensively in the literature.
For more details,  readers may be referred to \cite{He2013,FNP2,Nikiforov2014,Nikiforov2011,NY2,Yuan2014,Zhao2021}.
It is of interest to consider this problem for $F_{a_1, \cdots, a_k}$.

The \emph{adjacency matrix}
of $G$  is the $n\times n$ matrix $A(G)=(a_{ij})$, where
$a_{ij}=1$ if $v_i$ is adjacent to $v_j$, and $0$ otherwise.
Moreover, the matrix $Q(G)=D(G)+A(G)$ is known as  the \emph {signless Laplacian matrix} of $G$,  where $D(G)$ is the degree diagonal matrix of $G$.  The largest eigenvalues of $A(G)$ is called the \emph{spectral radius} of $G$.
The largest eigenvalues of $Q(G)$ is called the \emph{signless Laplacian spectral radius} of $G$ and denoted by $q(G)$.

In fact, if $k=1$ and $a_1=1$, Nikiforov \cite{Nikiforov2007} determined the maximum spectral radius of all $F_{1}$-free graphs of order $n$  and proved that  the only balance bipartite graphs has the maximum spectral radius; while  He, Jin and Zhang \cite{He2013}  obtained a sharp bound for the signless Laplacian spectral radius of all $F_{1}$-free graphs of order $n$ and proved that the corresponding  extremal graphs are all complete bipartite graphs $K_{r,s}$ with $r+s=n$.   Later, if $k=1$ and $a_1\ge 2$, Nikiforov \cite{Nikiforov2010} and Yuan \cite{Yuan2014} determined the maximum spectral radius  and the signless spectral radius among $F_{a_1}$-free graphs of order $n$, respectively.
Recently,
 Cioab\u{a}, Feng, Tait and Zhang \cite{Cioaba2020} studied the Brualdi-Solheid-Tur\'{a}n
type problem for $F^{(k)}$-free graphs and determined the corresponding spectral extremal graphs, which can be viewed as a spectral analogue of Theorem~\ref{thm1.1}.
\begin{Theorem}\label{thm1.3}\cite{Cioaba2020}
Let $G$ be an $F^{(k)}$-free graph of order $n$ with  $k \geq 2$. If $G$ has the maximum spectral radius, then
$$G \in Ex(n, F^{(k)})$$
for sufficiently large $n$.
\end{Theorem}
 Inspired by the work of Cioab\u{a}, Feng, Tait, and Zhang \cite{Cioaba2020}, Zhao, Huang, and Guo \cite{Zhao2021} focused on
the maximum signless Laplacian spectral radius of all $F^{(k)}$-free graphs of order $n$ and proved the following  result.
\begin{Theorem}\label{thm1.5}\cite{Zhao2021}
Let  $G$ be an $F^{(k)}$-free graph of order $n$. If $k\geq 2$ and $n \geq 3k^2 - k - 2$,  then
$$q(G)\leq q(S_{n,k})$$ with equality if and only if $G=S_{n,k}$.
\end{Theorem}

Recently,   Li and Peng \cite{Li2021} proved the spectral result of all $H_{k,s}$-free graphs of order $n$, which can be viewed as a spectral analogue of Theorem~\ref{thm1.2}.
\begin{Theorem}\label{thm1.4}\cite{Li2021}
Let   $k\ge 2$ and $1\le s\le k$.  If $G$ is  an $H_{k,s}$-free graph of order $n$ with the maximum spectral radius,  then
$$G \in Ex(n, H_{k,s})$$  for sufficiently large $n$.
\end{Theorem}

 Furthermore, they proposed the following conjecture on the signless Laplacian spectral radius of $F_{a_1, \cdots, a_k}$-free graphs of order $n$  with $k\ge 2$ and $a_1=\cdots=a_k\geq2$.
 \begin{Conjecture}\label{Cor1.6}\cite{Li2021}
 Let $G$ be  an $F_{a_1, \cdots, a_k}$-free  graph  of  order $n$. If   $k\ge 2$ and $a_1=\cdots=a_k=a\geq2$, then  $$q(G)\leq q(S_{n,ka})$$
  for sufficiently large $n$ with equality if and only if $G=S_{n,ka}$.
\end{Conjecture}

Inspired by the above known results and the conjecture on the signless Laplacian spectral radius of $F_{a_1, \cdots, a_k}$-free graphs with $k\ge 2$ and $a_1=\cdots=a_k\geq2$, we studied the maximum signless Laplacian spectral radius of all $F_{a_1, \cdots, a_k}$-free graphs.  Combined with known results, we present the main result of this paper as follows.

 \begin{Theorem}\label{m-thm}
 Let   $G$ be   an $F_{a_1,\dots,a_k}$-free  graph  of  order $n\geq 8t^2-12t+9$ with  $ t=\sum _{i=1}^k a_i$.

  (1) \cite{He2013} If $k=1$ and $a_1=1$, then  $q(G)\leq q(S_{n,t})$ with equality if and only if  $G$ is a complete bipartite graphs $K_{r,s}$ with $r+s=n$.

  (2) \cite{Yuan2014} If $k=1$, $a_1\ge2$, and  $n\geq 110 t^2$, then  $q(G)\leq q(S_{n,t})$  with equality if and only if $G=S_{n,t}$.

  (3) \cite{Zhao2021} If $k\ge 2$ and $a_1=\cdots=a_k=1$, then $q(G)\leq q(S_{n,t})$ with equality if and only if $G=S_{n, t}$.

  (4) If $k\ge 2$ and $a_1\ge 2$, then $q(G)\leq q(S_{n,t})$ with equality if and only if $G=S_{n, t}$.

  \end{Theorem}

{ \bf Remark.} It is worth mentioning that the extremal graphs in Theorem~\ref{m-thm}~(4) are not the
same as those of Theorems~\ref{thm1.2} and \ref{thm1.4}. The extremal graphs in Theorems~\ref{thm1.2} and \ref{thm1.4} are only depend  on the number of  intersecting triangles and the number of all intersecting odd cycles; while  the extremal graph in Theorem~\ref{m-thm}~(4) is not only depend on the number of intersecting odd cycles, but also the lengths of all intersecting odd cycles.
%The rest part of this paper, we always assume that $k\ge 2$ and $a_1\ge 2$.

The rest of this paper is organized as follows. In Section~2, some known lemmas are presented.
In Section~3, we give  the proof of Theorem~\ref{m-thm}.

\section{Some Lemmas}

\begin{Lemma}\label{thm2.1}\cite{Erdos1959}
Let $k\geq3$.  If  $G$ is $P_k$-free graph of order $n$, then $e(G)\leq \frac{(k-2)n}{2}$ with equality if and only if $G$ is union of disjoint copies of $K_{k-1}$.
\end{Lemma}

\begin{Lemma}\label{thm2.2}\cite{Erdos1959}
Let $k\geq2$.  If  $G$ is  graph of order $n$ with no cycle longer than $k$, then $e(G)\leq \frac{k(n-1)}{2}$ with equality if and only if $G=L_{r,k-1}$ with $n=r(k-1)+1$.
\end{Lemma}

\begin{Lemma}\label{lem2.6}\cite{Dirac1952}
 If  $G$ is a graph  with $\delta(G)\geq2$, then $G$ contains a cycle of length at least $\delta(G)$+1.
\end{Lemma}

 We also need the  stability result on the disjoint paths.

\begin{Lemma}\label{lem2.3}\cite{Chen2019-1}
Let $H=\bigcup_{i=1}^k P_{2a_i}$  with   $k\geq2$,
$a_1\geq \cdots \geq a_k\geq1$, and $ t=\sum_{i=1}^k a_i$. If $\delta (G)\geq t-1$ and  $G$ is  an $H$-free connected graph  of order $n\geq2t$.
 then  one of the following holds:\\
 (1) $G\subseteq S_{n,t-1}$;\\
  (2) $F=2P_{2a_1}$ and $G=L_{r,t-1}$, where $n=r(t-1)+1$.
\end{Lemma}

The structure of graphs and the spectral radius of graphs will be stated as follows.
\begin{Lemma}\label{lem2.4}\cite{Nikiforov2014}
Let $ t\geq2$ and  $n> 5t^2$. Then \\
(1)  $q(S_{n,t}) >n+2t-2-\frac{2(t^2-t)}{n+2t-3}>n+2t-3$.\\
  (2) If $G$ is a graph of order $n$ with  $q(G)\geq q(S_{n,t})$, then  $e(G)\geq tn-t^2+1$.
\end{Lemma}

\begin{Lemma}\label{lem2.5}
If $n=r(t-1)+2$ with $r\geq1$ and $ t\geq3$,  then
$q(K_1\nabla L_{r,t-1})<n+2t-3$.
\end{Lemma}

\begin{Proof}
Let $G=K_1\nabla L_{r,t-1}$.  Note that the nonincreasing degree sequence $(d_1,\dots,d_n)$ of $G$ is $(n-1,n-1,t,\dots,t)$.
By \cite[Theorem~4.2]{Duan2013}, we have
\begin{eqnarray*}
    q(G) &=& \frac{d_1+2d_3-1+\sqrt{(2d_3-d_1+1)^2+8\sum_{i=1}^2(d_i-d_3)}}{2} \\
   % &=& \frac{n-1+2t-1+\sqrt{(2t-(n-1)+1)^2+16(n-1-t)}}{2}\\
    &=& \frac{n+2t-2+\sqrt{(n-2t-2)^2+16(n-2t-2)+16t+16}}{2}\\
  & \leq&\frac{n+2t-2+(n-2t-2)+\frac{3t+7}{2}}{2}\\
  &<& %n+\frac{3t}{4}-\frac{1}{4}
  n+2t-3.
  \end{eqnarray*}
  This completes the proof.
\end{Proof}

\section{Proof of Theorem~\ref{m-thm}}
 In order to prove Theorem~\ref{m-thm}, we only to prove the assertion  holds for  $F_{a_1, \cdots, a_k}$  with
$k\ge 2$ and $a_1\ge 2$. Hence we always assume that $F_{a_1, \cdots, a_k}$  is the graph with $k\ge 2$ and $a_1\ge 2$ in this section. Firstly, we present several technical lemmas.

%The ideas of the proof in Lemma~\ref{Lem3.1} and Lemma~\ref{Lem3.2} come from Proposition 8 and Lemma 9 in
%\cite{NY}, respectively.

\begin{Lemma}\label{lem3.1}
 Let $G$ be  an $F_{a_1,\dots,a_k}$-free  graph  of  order $n\geq 8t^2-12t+9 $ with $ t=\sum _{i=1}^k a_i$. If    $q(G)\geq q(S_{n,t})$,
then $\Delta(G)=n-1$.
\end{Lemma}

\begin{Proof} By the upper bound for the signless Laplacian spectral radius from Merris \cite{Merris1998},  there exists a  vertex  $u$ such that
%  Note that
%$$q(G)\leq \max _{v\in V(G)}\Big\{d(v)+\frac{1}{d(v)}\sum_{z\in N(v)} d(z)\Big\},$$
%which dates back to .
%\leq \frac{2e(G)}{n-1}+n-2.$$,
%Let $u$ be a vertex such that
\begin{eqnarray*}%\label{merris}
q(G) \leq \max _{v\in V(G)}\Big\{d(v)+\frac{1}{d(v)}\sum_{z\in N(v)} d(z)\Big\}=d(u)+\frac{1}{d(u)}\sum_{z\in N(u)} d(z).
\end{eqnarray*}
 Let $A=N(u)$ and $B=V(G)\backslash (N(u)\bigcup \{u\})$.
 Then $|A|+|B|+1=n$ and
 \begin{equation}\label{E1}
\begin{aligned}
 q(S_{n,t})\leq q(G) \leq d(u)+\frac{1}{d(u)}\sum_{z\in N(u)} d(z)=|A|+1+\frac{2e(A)+e(A,B)}{|A|}.
     \end{aligned}
 \end{equation}
  Next we show that $d(u)=n-1$.   Assume for a contradiction that $d(u)<n-1$. We  prove the following claims.

 {\bf Claim~1.} No vertex in $B$ is adjacent to every vertex in $A$.

Suppose that there exists a vertex $v\in B$ which is adjacent to every vertex in $A$. Since $G$ is   $F_{a_1,\dots,a_k}$-free, we have
  $G[A]$ is $P_{2a_1-1}\bigcup(\bigcup_{i=2}^k P_{2a_i})$-free. Otherwise we can add the path  $P_{2a_1-1}$ in $G[A]$ to a path of order $2a_1$ by deleting the first edge $e$ of $P_{2a_1-1}$  with $a_1\ge 2$ and adding two edges from $u$ to the  endvertices of  $e$ and thus we obtain an $F_{a_1,\dots,a_k}$, which is a contradiction. Hence $G[A]$ is $P_{2t-1}$-free, where $t=\sum_{i=1}^k{a_i}.$
  By Lemma~\ref{thm2.1},
    we have $2e(A)\leq (2t-3)|A|$.
 By (\ref{E1}) and $e(A,B)\leq |A||B|$, we have
 \begin{eqnarray*}
  q(S_{n,t}) &\leq&  |A|+1+\frac{2e(A)+e(A,B)}{|A|} \\
    &\leq& |A|+1+\frac{(2t-3)|A|+|A||B|}{|A|} \\
   &=& %|A|+1+|B|+2t-3 =
   n+2t-3<q(S_{n,t}),
 \end{eqnarray*}
 where the last inequality is from Lemma~\ref{lem2.4}~(1) with $8t^2-12t+9\ge 5t^2$,
 which is a contraction. This proves Claim~1.

 {\bf Claim~2.}  $|B|\leq2t^2-2t$.

 By Claim~1, we have $e(A,B)\leq (|A|-1)|B|$. In addition, since $G$ is $F_{a_1,\dots,a_k}$-free, we have
  $G[A]$ is $(\bigcup_{i=1}^k P_{2a_i})$-free, which implies that $G[A]$ is $P_{2t}$-free. By Lemma~\ref{thm2.1}, we have
  $$2e(A)\leq (2t-2)|A|.$$
  Together with (\ref{E1}), we have
  \begin{eqnarray*}
  q(S_{n,t}) &\leq&  |A|+1+\frac{2e(A)+e(A,B)}{|A|} \\
    &\leq& |A|+1+\frac{(2t-2)|A|+(|A|-1)|B|}{|A|} \\
   %&=&|A|+1+|B|+2t-2-\frac{|B|}{n-1-|B|} \\
   &=& n+2t-2-\frac{|B|}{n-1-|B|}\\
   &<&n+2t-2-\frac{|B|}{n+2t-3}.
 \end{eqnarray*}
Hence  $|B|\leq2t^2-2t$ follows from Lemma~\ref{lem2.4}~(1). % we have $q(S_{n,t}) >n+2t-2-\frac{2(t^2-t)}{n+2t-3}$.
%Hence we get  .
This proves Claim~2.

Let $A'$ be the set of all vertices in $A$  which  are adjacent to every vertex in $B$.

 {\bf Claim~3.}  $|A'|\geq |A|-2t^2+2t$.

Since  $2e(A)\leq (2t-2)|A|$ and
$$e(A,B)\leq |A'||B|+(|A|-|A'|)(|B|-1)=|A||B|-|A|+|A'|,$$  we have

 \begin{eqnarray*}
  q(S_{n,t}) &\leq&  |A|+1+\frac{2e(A)+e(A,B)}{|A|} \\
    &\leq& |A|+1+\frac{(2t-2)|A|+|A||B|-|A|+|A'|}{|A|} \\
   &=&n+2t-2-\frac{|A|-|A'|}{|A|}\\
   &\le & n+2t-2-\frac{|A|-|A^{\prime}|}{n+2t-3}.
 \end{eqnarray*}
Hence $|A'|\geq |A|-2t^2+2t$  follows from Lemma~\ref{lem2.4}~(1). %, we have $q(S_{n,t}) >n+2t-2-\frac{2(t^2-t)}{n+2t-3}$.
%Hence we get $|A'|\geq |A|-2t^2+2t$ after some algebra.
This proves Claim~3.

Let $G_1$ be the union of all components of $G[A]$ each of which contains at least a vertex in  $A'$,
and let $G_2$ be the union of the remaining components of $G[A]$. Set $n_1=|V(G_1)|$ and $n_2=|V(G_2)|$. Note that $G_2$ is also $(\bigcup_{i=1}^k P_{2a_i})$-free which implies that  $G_2$ is $P_{2t}$-free. By Lemma~\ref{thm2.1} again,
$e(G_2) \leq (t - 1)n_2$.

{\bf Claim~4.}  $G_1$ contains no cycle longer than $2t -3$.

Assume for a contradiction that $G$ contains a cycle $C_p$  with $p \geq 2t - 2$. Since  $G$ is $F_{a_1,\dots,a_k}$-free graph,
  $G[A]$ must be  $\bigcup_{i=1}^k P_{2a_i}$-free,   which implies  $p\le 2t-1$.
  %that  $G$ contains $F_{a_1,\dots,a_k}$  as a subgraph.
  Hence $ 2t - 2\leq p\leq 2t-1$.
By the definition of  $G_1$,  there must be a vertex $z\in A'$ that either belongs to $C_p$
or can be joined to a vertex of $C_p$ by a path $Q$ contained in $G_1$.
By Claims~2-3,
$$|A'|\geq |A| - 2t^2+2t = n -1 - |B| -2t^2+2t \geq n  -4t^2 +4t-1 \geq2t >2.$$  We  consider the following two cases.

{\bf Case~1.} $z\in V(C_p)$.

 There exists  a vertex $v\in A'\backslash V(C_p)$ by $|A^{\prime}|\geq2t> p$. Choose a vertex $w\in B$ and construct a path $P$ of order $2t$  whose first three successive
  vertices are $v, w, z $ respectively and all vertices of $P$ except $w$ are  in $A$.
   As a result, $G[A\bigcup \{w\}]$ contains $\bigcup_{i=1}^k P_{2a_i}$ as a subgraph, where $\bigcup_{i=1}^k V(P_{2a_i})=V(P)$ and the  first three successive
  vertices of $P_{2a_1}$ are  $v, w, z $ respectively. Thus $G$ contains $F_{a_1,\dots,a_k}$  as a subgraph, which is a contradiction.

 {\bf Case~2.} $z$ is joined to a vertex of $C_p$ by a path $Q$ contained in $G_1$.

   Choose a vertex $w\in B$.  If there exists a vertex $v\in A'\backslash (V(C_p)\bigcup V(Q))$, then we can get a path $P$ of order $2t$  whose first three successive
  vertices are $v, w, z $  respectively and  all vertices $P$ except $w$ are  in $A$.  Otherwise, we  can also construct a path $P$ of order $2t$  whose first three successive
  vertices are $z, w, v $  respectively and all vertices of $P$ except $w$ are in $A$. Then $G[A\bigcup \{w\}]$ contains $\bigcup_{i=1}^k P_{2a_i}$ as a subgraph, implying that  $G$ contains $F_{a_1,\dots,a_k}$  as a subgraph.  This is a contradiction. Hence Claim 4 holds.

  By Claim~4  and Lemma~\ref{thm2.2},  we have   $$2e(G_1)\leq (2t-3)n_1.$$
  Then
  \begin{eqnarray*}
    2e(A) &=& 2e(G_1)+2e(G_2) \\
     &\leq& (2t-3)n_1 +(2t-2)n_2=(2t-3)(n_1+n_2)+n_2\\
    &\leq& (2t-3)|A|+|A|-|A'|\leq(2t-3)|A|+2t^2-2t.
  \end{eqnarray*}
By (\ref{E1}),
   \begin{eqnarray*}
 % n+2t-2-\frac{2(t^2-t)}{n+2t-3}<
 q(S_{n,t}) &\leq&  |A|+1+\frac{2e(A)+e(A,B)}{|A|} \\
    &\leq& |A|+1+\frac{(2t-3)|A|+2t^2-2t+|A||B|}{|A|} \\
   &=&n+2t-3+\frac{2t^2-2t}{|A|}.%\\
   %&\le & n+2t-2-\frac{|A|-(2t^2-2t)}{n}.
 \end{eqnarray*}
 Hence it follow from Lemma~\ref{lem2.4} (1) that
 $$n+2t-2-\frac{2t^2-2t}{|A|}\le n+2t-2-\frac{2t^2-2t}{n+2t-3}\le n+2t-2 -1+\frac{2t^2-2t}{|A|}
$$
 which implies that $|A|\le 4t^2-4t$. Hence by Claim 3,
 $$n=1+|A|+|B|\leq 1+4t^2-4t+2t^2-2t< 8t^2-12t+9,$$
  which is a contradiction.
% $$1-\frac{2t(t-1)}{n+2t-3}< \frac{2t^2-2t}{|A|}\leq \frac{2t^2-2t}{n-2t^2+2t}.$$
%This is a contradiction for $n\geq 6t^2-6t$.
Hence $d(u)=n-1$ and it completes the proof.
\end{Proof}

\begin{Lemma}\label{lem3.2}
 Let $G$ be a $\bigcup_{i=1}^ k P_{2a_i}$-free  graph  of  order $n>t^2-t-1$ with  $k\geq2$,  $a_1\geq \cdots \geq a_k\geq1$, and $ t=\sum _{i=1}^k a_i$. If $e(G)\geq  (t-1)n-(t^2-t-1)$ , then there exists an induced subgraph $H$ of $ G$ such that $|V(H)| \geq n-(t^2-t-1)$, $\delta(H)\geq t-1$, and $d_H(v)\leq t-2$ for every vertex $v\in V(G)\backslash V(H)$.
    %, and $ t=\sum _{i=1}^k a_i$.  If $G$ is  a $\bigcup_{i=1}^ k P_{2a_i}$-free  graph  of  order $n$ and $e(G)\geq  (t-1)n-(t^2-t-1)$,
 \end{Lemma}

\begin{Proof}
 If $\delta(G)\geq t-1$, then let $H=G$. Otherwise we delete a vertex  $v_1$ with the minimum degree $\delta(G)\le t-2$ and obtain $G_1=G-v_1$. If
%Next we  assume that $\delta(G)\leq t-2$
$\delta(G_1)\ge t-1$, let $H=G_1$. Otherwise we may continue this way until $\delta(G_{r-1})\le t-2$ and $G_r$ obtained from $G_{r-1}$ by deleting a vertex with the minimum degree and $\delta (G_r)\ge t-1$. We  show that $r\le t^2-t-1$.

% We construct a  sequence of graphs $G=G_0\supseteq G_1\supseteq \dots\supseteq G_r$
%such that  $\delta(G_i)\leq t-2$ for  $0\leq i\leq r-1$  and $\delta(G_r)\geq t-1$  by deleting a vertex with the minimum degree.  Now
Since $G_r$ is $\bigcup_{i=1}^ k P_{2a_i}$-free, we have $G_r$ is $P_{2t}$-free.  By $|V(G_r)|=n-r$ and Lemma~\ref{thm2.1},
$$(t-1)(n-r)\geq e(G_r)\geq e(G)-r(t-2)\geq (t-1)n-(t^2-t-1)-r(t-2),$$
which implies  $r\leq t^2-t-1$.
Hence $|V(H)|\ge n-(t^2-t-1)$.
%Let $H=G_r$, where $r-1$ is the maximum value of $i$ such that $\delta(G_i)\leq t-2$. This implies that $\delta(H)\geq t-1$.
%From the procedure constructing $\{G_i\}^{r}_{i=0}$,  $d_H(v)\leq t-2$ for  every vertex $v\in V(G)\backslash V(H)$.
%
 %Suppose that we have reached such a graph $G_i$ for some $i$ and $\delta(G_i)\leq t-2$.
 %Choose $v\in V(G_i)$ with $d_{G_i}(v)=\delta(G_i)$ and let $G_{i+1}=G_i-v$.
%Since $G_i$ is $\bigcup_{i=1}^ k P_{2a_i}$-free, we have $G_i$ is $P_{2t}$-free. Note that $|V(G_i)|=n-i$. By Lemma~\ref{thm2.1},
%$$(t-1)(n-i)\geq e(G_i)\geq e(G)-i(t-2)\geq (t-1)n-(t^2-t-1)-i(t-2).$$
%Thus $i\leq t^2-t-1$. Let $H=G_r$, where $r-1$ is the maximum value of $i$ such that $\delta(G_i)\leq t-2$. This implies that $\delta(H)\geq t-1$.
From the procedure constructing $\{G_i\}^{r}_{i=1}$,  $d_H(v)\leq t-2$ for  every vertex $v\in V(G)\backslash V(H)$.
\end{Proof}

%\begin{Lemma}\label{lem3.3}
%Let $H$ be the induced subgraph of $G$ in Lemma~3.2  such that $H=\bigcup_{i=1}^{s}H_i$ and $|V(H_i)|\geq 2t+2$ for $i=1,\ldots,s$, where $H_1,\ldots,H_s$ are the components of $H$. Then $H_i\neq L_{r_i,t-1}$  for $i=1,\ldots,s$, where $|V(H_{i})|=r_i(t-1)+1$.
%\end{Lemma}
%
%\begin{Proof} Suppose there exist an $H_m$ such that $H_m= L_{r_m,t-1}$ with $V(H_{m})=r_m(t-1)+1$, then $e(H_m)=\frac{kh_m-1}{2}$. Since $G$ is $\bigcup_{i=1}^{k}P_{2a_i}$-free, we have $H_m$ and $H\setminus V(H_m)$ also $\bigcup_{i=1}^{k}P_{2a_i}$-free. Then $H_m$ and $H\setminus V(H_m)$ are also $P_{2t}$-free.
%By Lemma~\ref{thm2.1},
%\begin{equation}\label{E2}
% \begin{aligned}
%    e(H) = &e (H\setminus V(H_i))+e(H_i) \\
%    \leq&(t-1)(|V(H)|-|V(H_i)|)+\frac{t(|V(H_i)|+1)}{2}\\
%     =&(t-1)|V(H)|-\frac{(t-2)|V(H_i)|}{2}-\frac{t}{2}.
%  \end{aligned}
% \end{equation}
%On the other hand, by Lemma~\ref{lem3.2},
%\begin{equation}\label{E3}
%\begin{aligned}
%    e(H)& \geq e(G)-(n-|V(H)|)(t-2) \\
%    & =(t-1)n-(t^2-t-1)-(n-|V(H)|)(t-2) \\
%    &=n+(t-2)|V(H)|-(t^2-t-1).
%    \end{aligned}
% \end{equation}
%By (\ref{E2}) and (\ref{E3}), we get
%$$n-|V(H)|+\frac{(t-2)|V(H)|}{2}\leq\frac{2t^2-3t+2}{2}=\frac{(2t+1)(t-2)}{2}$$
%which implies
%$$(t-2)|V(H_i)|\leq(2t+1)(t-2),$$
%i.e., $$|V(H_i)|\leq 2t+1,$$
%which is a contradiction.
%\end{Proof}

\begin{Lemma}\label{lem3.4}
 Let $G_1$ be a $\bigcup_{i=1}^k P_{2a_i}$-free graph of order $n_1$  with $t=\sum_{i=1}^k a_i\ge 3$ and $G_2$ be a graph of order $n_2$ with $1\le n_2\le t^2+t-2$. If
   $G=K_1\nabla (G_1\bigcup G_2)$  is a graph of order $n$ with $n=n_1+n_2+1\ge 8t^2-12t+9$, then
 %  $t\geq 3$ and $G=K_1\nabla (G_1\bigcup G_2)$  be a graph of order $n\geq 8t^2-12t+7$,   where $V(G_1)\bigcap V(G_2)=\emptyset$ and $1\leq |V(G_2|\leq t^2+t-2$. If $G_1$ is $\bigcup_{i=1}^k P_{2a_i}$-free, then
  $q(G) < q(S_{n,t})$.
\end{Lemma}
\begin{Proof}
Without loss of generality, we assume that $G_2=K_{n_2}$.
%Let $u$ be the dominating vertex of $G$ and $|V(G_2)|=p$. Then $1\leq p\leq t^2+t-2$. Suppose that $q(G)\geq q(S_{n,t})$. We may  assume that $G_2=K_{p}$, otherwise, adding edges to $G$ to get such a graph whose signless Laplacian spectral radius will not decrease.
Since $G_1$ is $\bigcup_{i=1}^k P_{2a_i}$-free, we have $G_1$ is $P_{2t}$-free. By Lemma~\ref{thm2.1},
$$e(G_1)\leq (t-1)(n-n_2-1)$$
Thus,
$$e(G-V(G_2))=e(G_1)+n-n_2-1\leq t(n-n_2-1).$$
By \cite[Theorem~3.1]{Das2004},
\begin{eqnarray*}
 q(G-V(G_2)) &\leq & \frac{2e(G-V(G_2))}{|V(G-V(G_2))|-1}+|V(G-V(G_2))|-2\\
   &\leq & \frac{2t(n-n_2-1)}{n-n_2-1}+n-n_2-2 \\
  &\leq&  n+2t-n_2-2.
\end{eqnarray*}
Let ${\bf x}=(x_1,\ldots, x_n)^T$ be a unit positive eigenvector of $Q(G)$ corresponding  to the signless Laplacian spectral radius $q(G)$. By symmetry, all components of $x$ corresponding to the vertices of $G_2$  are the same, say $a$. By the eigenequations of $Q(G)$ on $u$ with the maximum degree and a vertex of $G_2$, respectively, we have

 \begin{eqnarray*}
   (q(G)-n+1)x_u&=&\sum_{ v\in V(G)\setminus \{u\}}x_v \leq \sqrt{(n-1)\sum_{ v\in V(G)\setminus \{u\}}x^2_{v}} =\sqrt{(n-1)(1-x^2_{u})},\\
  (q(G)-n_2)a&=&(n_2-1)a+x_{u}.
     \end{eqnarray*}
Note that $q(G)\geq q(S_{n,t})> n+2t-3$ and $n_2\leq t^2+t-2$,  we have

 \begin{eqnarray*}
  x^2_{u}\leq \frac{n-1}{(q(G)-n+1)^2+n-1}\leq \frac{n-1}{n-1+4(t-1)^2}< 1-\frac{4(t-1)^2}{n+4t^2},
     \end{eqnarray*}
     and \begin{eqnarray*}
       a =\frac{x_u}{q(G)-2n_2+1}\leq \frac{x_u}{n+2t-2n_2-2}\leq \frac{x_u}{n-2t^2+2}.
       \end{eqnarray*}
Then
\begin{eqnarray*}
q(G)&=&\sum_{ij\in E(G)} (x_i+x_j)^2=\sum_{ij\in E(G)\setminus E(G_2)}(x_i+x_j)^2+n_2(a+x_u)^2+\sum_{ij\in E(G_2)}(x_i+x_j)^2\\
&<&q(G-V(G_2))+n_2(a+x_u)^2+2n_2(n_2-1)a^2\\
&\leq&n+2t-n_2-2+n_2\left(1+\frac{1}{(n-2t^2+2)^2}+\frac{2}{n-2t^2+2}+\frac{2(n_2-1)}{(n-2t^2+2)^2}\right)x^2_u\\
&<&n+2t-n_2-2+n_2\left(1+\frac{3}{n-2t^2+2} \right)\left(1-\frac{4(t-1)^2}{n+4t^2}\right)\\
&<&n+2t-2-\left( \frac{4(t-1)^2}{n+4t^2}-\frac{3}{n-2t^2+2}\right)\\
&<&n+2t-2-\frac{2t(t-1)}{n-2t^2+2}\\
&<&n+2t-2-\frac{2t(t-1)}{n+2t-3}<q(S_{n,t}).
      \end{eqnarray*}
This completes the proof.
\end{Proof}

\vspace{2mm}
Now we are ready to prove Theorem~\ref{m-thm}.

\vspace{2mm}
\noindent {\bf Proof of Theorem~\ref{m-thm}.}
Suppose that $q(G)\geq q(S_{n,t})$. We will show that $G=S_{n,t}$. By Lemma~\ref{lem3.1}, $\Delta(G)=n-1$. Let $u\in V(G)$ be the vertex with the maximum degree $\Delta(G)$, i.e., $d(u)=n-1$.  By Lemma~\ref{lem2.4}~(2), $e(G)\geq tn-t^2+1$, which implies that
\begin{eqnarray*}
  e(G-u)&=& e(G)-n+1\\
  &\geq& (t-1)n-t^2+2\\
  &=&(t-1)(n-1)-(t^2-t-1).
\end{eqnarray*}
Note that $G-u$ is $\bigcup_{i=1}^k P_{2a_i}$-free, otherwise $G$ contains $F_{a_1,\dots,a_k}$as a subgraph. By Lemma~\ref{lem3.2} , there exists an induced subgraph $H$ of $G-u$ such that $\delta(H)\geq t-1$, $|V(H)|=n_1 \geq n-(t^2-t)$, and $d_H(v)\leq t-2$ for every vertex $v\in V(G)\backslash (V(H)\bigcup \{u\})$. Let $H=\bigcup_{i=1}^s H_i$ and $|V(H_i)|=h_i$, where $H_i$ is a component of $H$ for $i=1, \ldots, s.$ %For $i=1, \ldots, s$, since $\delta(H_i)\geq \delta(H)\geq t-1\geq 2.$ By Lemma~\ref{lem2.6},  $H_i$ contains a cycle of length at least $t$.

{\bf Claim 1.} Every component of $G-u$ contains at most one graph of ${H_1, \ldots, H_s}$ as an induced subgraph.

Note that $\delta(H_i)\geq t-1\geq2$ for $i=1,\dots,s$. By Lemma~\ref{lem2.6}, $H_i$ contains a cycle of length at least $t$  for $i=1,\dots,s$.  In fact, if there is a component of $G-u$ containing at least two graphs $H_i$ and $H_j$ as an induced subgraph for $1\leq i\neq j \leq s$, then $G-u$ contains $P_{2t+1}$ as a subgraph and thus  $G$ contains $F_{a_1,\ldots,a_k}$ as a subgraph, which is a contradiction. This proves Claim~1.

By Claim 1, we  let $T_i$ be the component of $G-u$ containing $H_i$ as an induced subgraph, and $G-u=(\bigcup_{i=1}^{s}T_i)\bigcup T_0,$ where $T_0$ is the union of the remaining components of $G-u$. Note that $G-u$ is $\bigcup_{i=1}^k P_{2a_i}$-free, we have $T_i$ is $P_{2t}$-free  for $i=0,\ldots, s$.

{\bf Claim 2.} $T_0= \emptyset$.

In fact, if $T_0 \neq \emptyset$, then
\begin{eqnarray*}
    1 \leq|V(T_0)|&=&n-1-\sum_{i=1}^{s}|V(T_i)| \leq n-1-\sum_{i=1}^{s}|V(H_i)| \\
    & =&n-1-n_1\leq t^2-t-1.
  \end{eqnarray*}
By Lemma~\ref{lem3.4}, $q(G)<q(S_{n,t})$, which is a contradiction. This proves Claim~2.

%By Claim 2, we have $$G-u=\bigcup_{i=1}^{t}T_i.$$

{\bf Claim 3.} $h_i\geq 2t$ for $i=1,\ldots, s$.

Suppose  there exist an $h_i$ such that $h_i\leq 2t-1$.  Since $\delta(H_i)\geq \delta(H)\geq t-1$, we have $h_i\geq t$. Then
\begin{eqnarray*}
    t \leq h_i \leq  |V(T_i)|\leq h_i+|V(G-u)\setminus V(H)| \leq 2t-1+t^2-t-1=t^2+t-2.
  \end{eqnarray*}
By Lemma~\ref{lem3.4}, $q(G)< q(S_{n,t})$, which is a contradiction. This proves Claim~3.

By the definition of $L_{r,t}$, we have the following claim directly.

{\bf Claim~4.}  For any fixed $1\leq i\leq s$, if $H_i=L_{r_i,t-1}$ with $h_i=r_i(t-1)+1$, then $H_i$ contains $P_{2t-1}$ as a subgraph.

{\bf Claim~5.} For any fixed $1\leq i\leq s$, if $H_i$ is a subgraph of $S_{h_i,t-1}$, then $H_i$ contains $P_{2t-1}$ as a subgraph. Moreover,  $T_i$ is subgraph of $S_{|V(T_i)|,t-1}$.

If $H_i$ is a subgraph of $S_{h_i,t-1}$, then there exists $I_i\subseteq V(H_i)$ of size $h_i-t+1$ such that $I_i$ induces an independent set of $H_i$. Since $\delta(H_i)\geq t-1$, every vertex in $I_i$ is adjacent to every vertex in $V(H_i)\backslash I_i$.  Then $H_i$ contains a path $P$ of order $2t-1$ with both endvertices in $I_i$ as a subgraph.
If $(V(T_i)\backslash V(H_i))\bigcup I_i$ induces at least an edge, then we can get a path of order $2t$ from $P$. Hence $G$ contains $F_{a_1,\dots,a_k}$ as a subgraph, which is  a contradiction.
This implies that $(V( T_i)\backslash V(H_i))\bigcup I_i$ is an independent set and thus   $T_i$ is subgraph of $S_{|V(T_i)|,t-1}$. This proves Claim~5.

Note that $q(G) \geq q(S_{n,t})$. Since $H_i$ is $\bigcup_{i=1}^k P_{2a_i}$-free and $\delta(H_i)\geq \delta(H)\geq t-1$, by Lemmas~\ref{lem2.3},
$H_i$ is a subgraph of $S_{h_i,t-1}$ or $H_i=L_{r_i,t-1}$ with $h_i=r_i(t-1)+1$ and $k=2, a_1=a_2$ for $i=1,\dots,s$. If $s\geq2$, then by Claims~3 and 4, $G-u$ contains $2P_{2t-1}$ as a subgraph and thus $G$ contains $F_{a_1,\dots,a_k}$ as a subgraph, which is a contradiction. So $s=1$. This implies that
$$G-u=T_1,$$
and $H_1$ is an induced graph of $T_1$, where $H_1$  is a subgraph of $S_{h_1,t-1}$ or $H_1=L_{r_1,t-1}$ with $h_1=r_1(t-1)+1$  and $k=2, a_1=a_2$.

First suppose that  $H_1$  is a subgraph of $S_{h_1,t-1}$. By Claim~5, $T_1$ is  a subgraph of $S_{|V(T_1)|,t-1}$, i.e., $G$ is a subgraph of $S_{n,t-1}$. If $G$ is a proper subgraph of $S_{n,t}$, then it follows from Perron-Fronbenius theorem that
$q(G)<q(S_{n,t})$, which is also a contradiction. Hence $G=S_{n,t}.$

Next suppose that  $H_1=L_{r_1,t-1}$ with $h_1=r_1(t-1)+1$  and $k=2, a_1=a_2$. If $H_1=T_1$, then $G=K_1\nabla L_{r_1,t-1}$ with $n=r_1(t-1)+2$. By Lemma~\ref{lem2.5},
$q(G)<n+2t-3<q(S_{n,t})$, which is a contradiction. Thus $H_1$ is a proper subgraph of $T_1$.
Let $H'=T_1-V(H_1)$ and $|V(H')|=n_2$.
Since $d_{H_1}(v)\leq t-2$ for every vertex $v\in V(H')$,
we have$$e(V(H'),V(H_1))\leq (t-2) n_2\leq(t-2)(t^2-t-1).$$
%Then
%\begin{eqnarray*}
%% \nonumber to remove numbering (before each equation)
%  e(H) &\leq &e(L_{t_1,t_2,h,h+1})
%    =\frac{t_1h^2+t_2(h+1)^2+n_1-1}{2} \\
%   &\leq& \frac{t_1h(h+1)+t_2(h+1)^2+n_1-1}{2}
%   = \frac{(h+2)(n_1-1)}{2}.
%\end{eqnarray*}
%Furthermore, since  $q(G)\geq q(S_{n,h})$, by Lemma~\ref{Lem2.2}~($iii$), we have $e(G)>h(n-h)$.
Then
 \begin{eqnarray*}
     %\nonumber to remove numbering (before each equation)
        e(H')&=& e(T_1)-e(H_1)-e(V(H'),V(H_1))  \\
       & > & (t-1)n-t^2+2-\frac{t(n_1-1)}{2}-(t-2)(t^2-t-1)\\
       &=&(t-1)n-t^2+2-\frac{t(n-n_2-2)}{2}-(t-2)(t^2-t-1)\\
        & = &\frac{(t-2)n-(2t^3-4t^2)}{2}+\frac{t n_2}{2}
         > \frac{t n_2}{2}.
  \end{eqnarray*}
By Lemma~\ref{thm2.1},  $H'$ contains $P_{t+2}$  as a subgraph. Together with Claim~4, $T_1$ contains $P_{2t-1}\bigcup P_{t+2}$ as a subgraph. This implies that $G$ contains
 $F_{a_1,a_2}$ as a subgraph, which is a contradiction. This completes the proof.
 \QEDB

\vspace{3mm}

%\noindent{\bf Acknowledgements.} The authors would like to thank the anonymous referee for many helpful and constructive suggestions to an earlier version of this paper, which results in a great improvement.


\begin{thebibliography}{12}

%\bibitem{CDS}
%D.M. Cvetkovi\'{c}, M. Doob, H. Sachs, Spectra of Graphs, %Theory and Application,
%third ed., Johann Ambrosius Barth, Heidelberg, 1995.
%\bibitem{CRS97} D.M. Cvetkovi\'c, P. Rowlinson, S.K. Simi\'c,  Eigenspaces of Graphs, Cambridge Univ.
%Press, Cambridge, 1997.
%\bibitem{Z}
%B. Zhou, Signless Laplacian spectral radius and Hamiltonicity, Linear Algebra Appl. 432 (2010) 566--570.



%\bibitem{BLL}   B. Bollob\'{a}s, J. Lee, S. Letzter, Eigenvalues of subgraphs of the cube, preprint, arXiv:1605.06360, 2016.

%\bibitem{BK}   N. Bushaw and N. Kettle, Tur\'{a}m numbers of multiple paths and equibipartite forests, Combinatorics, Probability and
%Computing 20 (2011)  837--853.

%\bibitem{BZ}   A. Berman,  X.-D. Zhang, On the spectral radius of graphs with cut vertices, J. Comb. Theory Ser. B 83 (2001) 233--240.
% \bibitem{Ali1996} A.A. Ali, W. Staton, On extremal graphs with no long paths, Electron. J. Combin. 3 (1996), R20.

 \bibitem{BM}  J.A. Bondy,  U.S.R. Murty,   Graph Theory,  Springer, New York, 2007.

 % \bibitem{BK}   N. Bushaw, N. Kettle, Tur\'{a}n numbers of multiple paths and equibipartite forests, Combin. Probab. Comput. 20 (2011)  837--853.

%\bibitem{campos2015} V.~Campos, R.~Lopes, A proof for a conjecture of Gorgol, Electron. Notes  Discrete Math. 50 (2015) 367--372.

\bibitem{Chen2019-1}  M.-Z. Chen, X.-D. Zhang, Erd\H{o}s--Gallai stability theorem for linear forests,  Discrete Math. 342 (2019) 904--916.

\bibitem{Chen2019-2}  M.-Z. Chen, A.-M. Liu, X.-D. Zhang, Spectral extremal results with forbidding linear forests, Graphs Combin. 35 (2019) 335--351.

\bibitem{Chen2021} M.-Z. Chen, A-M. Liu, X.-D. Zhang. On the spectral radius of graphs without a star forest.
Discrete Math. 344 (2021), 112269.

 \bibitem{Cioaba2020} S. Cioab\v{a}, L. Feng, M. Tait, X.-D. Zhang, The maximum spectral radius of graphs
without friendship subgraphs, Electron. J. Combin. 27 (4) (2020) P4.22.

\bibitem{Das2004} K. Das, Maximizing the sum of the squares of the degrees of a graph, Discrete Math. 285 (2004) 57--66.

\bibitem{Duan2013}  X. Duan, B. Zhou, Sharp bounds on the spectral radius of a nonnegative
matrix, Linear Algebra Appl. 439 (2013) 2961--2970.

\bibitem{FNP2}  M.A.A. de Freitas, V. Nikiforov, L. Patuzzi, Maxima of the $Q$-index: graphs with no $K_{s,t}$, Linear Algebra Appl. 496 (2016) 381--391.


 \bibitem{Dirac1952} G.A. Dirac, Some theorems on abstract graphs, Proc. London Math. Soc. 2 (1952) 69--81.


% \bibitem{DPS}  D. Cvetkovi\'{c}, P. Rowlinson, S. Simi\'{c}, Signless Laplacian of finite graphs, Linear Algebra Appl. 423 (2007) 155--171.


 \bibitem{Erdos1959}  P. Erd\H{o}s,  T. Gallai, On maximal paths and circuits of graphs. Acta Math. Acad. Sci. Hungar. 10 (1959) 337--356.

\bibitem{Erdos1995}  P. Erd\H{o}s, Z. F\"{u}redi, R.J. Gould, D.S. Gunderson, Extremal graphs for intersecting triangles, J. Combin. Theory. Ser. B, 64 (1995) 89--100.
 %    \bibitem{FY}    L. Feng, G.  Yu, On three conjectures involving the signless Laplacian spectral radius of graphs, Publ. Inst. Math. (Beograd). 85 (2009) 35--38.

%  \bibitem{FNP}  M.A.A. de Freitas, V. Nikiforov, L. Patuzzi, Maxima of the $Q$-index: forbidden 4-cycle and 5-cycle, Electron. J. Linear Algebra 26
%      (2013) 905--916.



% \bibitem{Gorgol}  I. Gorgol, Tur\'{a}n number for disjoint copies of graphs, Graphs Combin. 27 (2011) 661-667.

  \bibitem{He2013}   B. He,   Y.-L. Jin,  X.-D. Zhang, Sharp bounds for the signless Laplacian spectral radius in terms of clique number, Linear Algebra Appl. 438 (2013)  3851--3861.

  \bibitem{Hou2018}  X. Hou, Y. Qiu, B. Liu, Tur\'{a}n number and decomposition number of intersecting odd cycles, Discrete Math. 341 (2018) 126--137.
%\bibitem{mantel1907}W. Mantel, Problem 28, soln. by H. Gouwentak, W. Mantel, J. Teixeira de Mattes,
%F. Schuh and W.A. Wythoff. Wiskundige Opgaue,i 10 (1907), 60--61.
   \bibitem{Li2021} Y. Li, Y. Peng, The spectral radius of graphs with no intersecting odd cycles, arXiv:2106.00587v1.
%
\bibitem{Merris1998}   R. Merris, A note on Laplacian graph eigenvalues, Linear Algebra Appl. 285 (1998) 33--35.

\bibitem{Nikiforov2007} V. Nikiforov, Bounds on graph eigenvalues II, Linear Algebra Appl. 427 (2007) 183--189.


\bibitem{Nikiforov2010}   V. Nikiforov, The spectral radius of graphs without paths and cycles of specified length, Linear Algebra Appl. 432 (2010) 2243--2256.

\bibitem{Nikiforov2010-2} V. Nikiforov, A contribution to the Zarankiewicz problem, Linear Algebra Appl. 432 (2010) 1405--1411.


%\bibitem{Nikiforov3}  V. Nikiforov, A Spectral Erd\H{o}s--Stone--Bollob\'{a}s Theorem, Combin. Probab. Comput. 18 (2009) 455--458.

\bibitem{Nikiforov2011}  V. Nikiforov, Some new results in extremal graph theory, Surveys in combinatorics 2011, 141--181, London Math. Soc. Lecture Note Ser., 392, Cambridge Univ. Press, Cambridge, 2011.

  \bibitem{Nikiforov2014} V. Nikiforov, X. Yuan,  Maxima of the $Q$-index: graphs without long paths, Electron. J. Linear Algebra 27 (2014) 504--514.

   \bibitem{NY2}  V. Nikiforov, X. Yuan, Maxima of the $Q$-index: forbidden even cycles, Linear Algebra Appl. 471 (2015) 636--653.
%
%   \bibitem{Wilf}   H.S. Wilf, Spectral bounds for the clique and independence numbers of graphs, J. Combin. Theory Ser. B 40 (1986) 113--117.

 %   \bibitem{RZ}  R. Xing,   B. Zhou, On the least eigenvalue of cacti with pendant vertices,  Linear Algebra Appl. 438 (2013) 2256--2273.
%
%  \bibitem{Stanley}  R.P. Stanley, A bound on the spectral radius of graphs with $e$ edges, Linear Algebra Appl. 87 (1987) 267--269.


  % \bibitem{Tait2017} M. Tait, J. Tobin, Three conjectures in extremal spectral graph theory, J. Combin. Theory Ser. B 126 (2017) 137--161.
%
%  \bibitem{Tait2019}  M. Tait, The Colin de Verdi\`{e}re parameter, excluded minors, and the spectral radius, J. Combin. Theory Ser. A 166 (2019) 42--58.
   %\bibitem{Stevanovic}   D. Stevanov\'{i}c, Spectral radius of graphs, Academic Press, Amsterdam (2015).

%  \bibitem{Turan}  P. Tur\'{a}n, On an extremal problem in graph theory, Mat. Fiz. Lapok 48 (137) (1941) 436--452.

    % \bibitem{TT}  M. Tait, J. Tobin, Three conjectures in extremal spectral graph theory, J. Combin. Theory Ser. B 126 (2017) 137--161.


  %   \bibitem{YGW} G. Yu, S.-G. Guo, Y. Wu,  Maxima of the $Q$-index for outer-planar graphs, Linear and Multilinear Algebra 63 (2015) 1837--1848.

\bibitem{Simonovits1968} M.~Simonovits, A method for solving extremal problems in graph theory, stability problems, Theory of Graphs (Proc. Colloq., Tihany, 1966) (1968) 279-319.

  \bibitem{Yuan2018}  L.-T. Yuan, Extremal graphs for the $k$-flower, J. Graph Theory 89 (2018) 26--39.

  \bibitem{Yuan2014}  X. Yuan, Maxima of the $Q$-index: forbidden odd cycles, Linear Algebra Appl. 458 (2014) 207--216.


%  \bibitem{YZ} L.-T.~Yuan, X.-D~Zhang, The Tur\'{a}n number of  disjoint copies of paths,  Discrete Math. 340 (2017) 132--139.

   \bibitem{Zhao2021} Y. Zhao, X. Huang, H. Guo, The signless Laplacian spectral radius of graphs with
no intersecting triangles, Linear Algebra Appl. 618 (2021) 12--21.

%[2] B. Bollob¨¢s, Modern Graph Theory, Grad. Texts in Math., vol. 184, Springer-Verlag, New York,
%1998.
%[3] D. Cvetkovi?, Spectral theory of graphs based on the signless Laplacian, Research Report, 2010,
%available at: http://www.mi.sanu.ac.rs/projects/signless_L_reportApr11.pdf.






\end{thebibliography}
\end{document}